\documentclass[a4paper,10 pt]{article}
\usepackage{amsmath}
\usepackage{amssymb}
\usepackage{amsthm}
\usepackage{amsfonts}
\usepackage{mathrsfs}
\usepackage{enumerate}
\begin{document}
\title{\Large $(2,3)$-GENERATION OF THE SPECIAL LINEAR GROUPS OF DIMENSIONS $9$, $10$ and $11$}
\author{\Large E. Gencheva, Ts. Genchev and K. Tabakov}
\date{}
\maketitle
\begin{abstract}
In the present paper we prove that the groups $PSL_{n}(q)$ ($n = 9$, $10$ or $11$) are $(2,3)$-generated for any $q$. Actually, we provide  explicit generators $x_{n}$ and $y_{n}$ of respective orders $2$ and $3$, for the  special linear group $SL_{n}(q)$.
\end{abstract}
\vspace{16pt}
\indent\indent\indent\textbf{1.\,\,Introduction.} $(2,3)$-generated groups are those groups which can be generated by an involution and an element of order $3$ or, equivalently, they appear to be homomorphic images of the famous modular group $PSL_{2}(\mathbb{Z})$. It is known that many series of finite simple groups are 
$(2,3)$-generated. Most powerful result of Liebeck-Shalev and L\"{u}beck-Malle states that all finite simple groups, except the infinite families $PSp_{4}(2^{m})$, $PSp_{4}(3^{m})$, $^{2}B_{2}(2^{2m+1})$, and a finite number of other groups, are $(2,3)$-generated (see \cite{12}). We have especially focused our attention to the projective special linear groups defined over finite fields. Many authors have been investigated the groups $PSL_{n}(q)$ with respect to that generation property. $(2,3)$-generation has been proved in the cases $n=2$, $q\neq 9$ \cite{7}, $n=3$, $q\neq 4$ \cite{5},\cite{2}, $n=4$, $q\neq 2$ \cite{14}, \cite{13}, \cite{8}, \cite {10}, $n=5$, any $q$ \cite{17}, \cite{9}, $n=6$, any $q$ \cite{16}, $n=7$, any $q$ \cite{15}, $n=8$, any $q$ \cite{6}, $n\geq 5$, odd $q\neq 9$ \cite{3}, \cite{4}, and $n\geq 13$, any $q$ \cite{11}. In this way the only cases that still remain open are those for $9\leq n \leq12$, even $q$ or $q=9$. In the present work we continue our investigation by considering the next portion of infinite series of finite special linear groups and their projective images. We shall treat the groups $SL_{9}(q)$ and $SL_{10}(q)$ simultaneously. Based on the results obtained below, we deduce the following.\\ 

\indent\indent\textbf{Theorem.}  \emph{The groups $SL_{9}(q)$, $SL_{10}(q)$, $SL_{11}(q)$ and their simple projective images $PSL_{9}(q)$, $PSL_{10}(q)$ and $PSL_{11}(q)$ are $(2,3)$-generated for all $q$ }.\\

\indent\indent\textbf{2.\,\,Proof of the Theorem.} First let $G = SL_{n}(q)$ and $\overline{G} = G/Z(G) = PSL_{n}(q)$, where $n = 9$ or $10$, and $q = p^{m}$  for a prime number $p$. Set $Q = q^{n-1}-1$ if $q \neq 3, 7$ and $Q = (q^{n-1}-1)/2$ if  $q = 3, 7$.  The group $G$ acts (naturally) on the left on the $n$-dimensional column vector space $V = F^{n}$ over the field $F = GF(q)$. We denote by $v_{1}$, . . . , $v_{n}$ the standard base of the space $V$, i.e. $v_{i}$ is a column which has $1$ as its $i$-th coordinate, while all other coordinates are zeros.\\
\indent\indent\ We shall need the following result, which can be easily obtained by the list of maximal subgroups of $G$ given in \cite{1} and simple arithmetic considerations using (for example) Zsigmondy's well-known theorem. \\

\indent\indent\textbf{Lemma 1.} \emph{For any maximal subgroup $M$ of the group $G$ either it stabilizes one-dimensional subspace or hyperplane of $V$ ($M$ is reducible on the space $V$) or $M$ has no element of order $Q$}.\\

\indent\indent \textbf{2.1.} We suppose first that $q \neq 2, 4$ if $n = 9$ and $q > 4$ if $n = 10$. Let choose an element $\omega$ of order $Q$ in the multiplicative group of the field $GF(q^{n-1})$ and set 
\begin{center}
$f_{n}(t) = (t - \omega)(t - \omega^{q})(t - \omega^{q^{2}})(t - \omega^{q^{3}}) . . . (t - \omega^{q^{n-3}}) (t - \omega^{q^{n-2}}) = t^{n-1} - \alpha_{1}t^{n-2} + \alpha_{2}t^{n-3} - \alpha_{3}t^{n-4} + . . . + (-1)^{n-2}\alpha_{n-2}t + (-1)^{n-1}\alpha_{n-1}$.
\end{center}
Then $f_{n}(t) \in F[t]$ and the polynomial $f_{n}(t)$ is irreducible over the field $F$. Note that  $\alpha_{n-1} = \omega^\frac{q^{n-1} - 1}{q - 1}$  has order  $q - 1$ if $q \neq 3, 7$, $\alpha_{n-1} = 1$ if $q = 3$, and  $\alpha_{n-1}^{3} = 1 \neq \alpha_{n-1}$  if  $q = 7$.\\
\indent\indent Now let
\begin{center}
\[ x_{9} = \left[ \begin{array}{ccccccccc}
-1 & 0 & 0 & 0 & 0 & 0 & \alpha_{5}\alpha_{8}^{-1} & 0 & \alpha_{5}\\
0 & -1 & 0 & 0 & 0 & 0 & \alpha_{4}\alpha_{8}^{-1} & 0 & \alpha_{4}\\
0 & 0 & 0 & -1 & 0 & 0 & \alpha_{3}\alpha_{8}^{-1} & 0 & \alpha_{6}\\
0 & 0 & -1 & 0 & 0 & 0 & \alpha_{6}\alpha_{8}^{-1} & 0 & \alpha_{3}\\
0 & 0 & 0 & 0 & -1 & 0 & \alpha_{2}\alpha_{8}^{-1} & 0 & \alpha_{2}\\
0 & 0 & 0 & 0 & 0 & 0 & \alpha_{1}\alpha_{8}^{-1} & -1 & \alpha_{7}\\
0 & 0 & 0 & 0 & 0 & 0 & 0 & 0 & \alpha_{8}\\
0 & 0 & 0 & 0 & 0 & -1 & \alpha_{7}\alpha_{8}^{-1} & 0 & \alpha_{1}\\
0 & 0 & 0 & 0 & 0 & 0 & \alpha_{8}^{-1} & 0 & 0\\
 \end{array} \right],\]
\[y_{9} = \left[ \begin{array}{ccccccccc}
0 & 0 & 1 & 0 & 0 & 0 & 0 & 0 & 0\\
1 & 0 & 0 & 0 & 0 & 0 & 0 & 0 & 0\\
0 & 1 & 0 & 0 & 0 & 0 & 0 & 0 & 0\\
0 & 0 & 0 & 0 & 0 & 1 & 0 & 0 & 0\\
0 & 0 & 0 & 1 & 0 & 0 & 0 & 0 & 0\\
0 & 0 & 0 & 0 & 1 & 0 & 0 & 0 & 0\\
0 & 0 & 0 & 0 & 0 & 0 & 0 & 0 & 1\\
0 & 0 & 0 & 0 & 0 & 0 & 1 & 0 & 0\\
0 & 0 & 0 & 0 & 0 & 0 & 0 & 1 & 0\\
 \end{array} \right];\]
\[ x_{10} = \left[ \begin{array}{cccccccccc}
0 & 0 & 0 & -1 & 0 & 0 & 0 & \alpha_{2}\alpha_{9}^{-1} & 0 & \alpha_{3}\\
0 & 0 & 0 & 0 & 0 & -1 & 0 & \alpha_{4}\alpha_{9}^{-1} & 0 & \alpha_{7}\\
0 & 0 & -1 & 0 & 0 & 0 & 0 & \alpha_{5}\alpha_{9}^{-1} & 0 & \alpha_{5}\\
-1 & 0 & 0 & 0 & 0 & 0 & 0 & \alpha_{3}\alpha_{9}^{-1} & 0 & \alpha_{2}\\
0 & 0 & 0 & 0 & 0 & 0 & 0 & \alpha_{1}\alpha_{9}^{-1} & -1 & \alpha_{8}\\
0 & -1 & 0 & 0 & 0 & 0 & 0 & \alpha_{7}\alpha_{9}^{-1} & 0 & \alpha_{4}\\
0 & 0 & 0 & 0 & 0 & 0 & -1 & \alpha_{6}\alpha_{9}^{-1} & 0 & \alpha_{6}\\
0 & 0 & 0 & 0 & 0 & 0 & 0 & 0 & 0 & \alpha_{9}\\
0 & 0 & 0 & 0 & -1 & 0 & 0 & \alpha_{8}\alpha_{9}^{-1} & 0 & \alpha_{1}\\
0 & 0 & 0 & 0 & 0 & 0 & 0 & \alpha_{9}^{-1} & 0 & 0\\
 \end{array} \right],\] 
\[y_{10} = \left[ \begin{array}{cccccccccc}
1 & 0 & 0 & 0 & 0 & 0 & 0 & 0 & 0 & 0\\
0 & 0 & 1 & 0 & 0 & 0 & 0 & 0 & 0 & 0\\
0 & 0 & 0 & 0 & 0 & 0 & 1 & 0 & 0 & 0\\
0 & 0 & 0 & 0 & 0 & 1 & 0 & 0 & 0 & 0\\
0 & 0 & 0 & 1 & 0 & 0 & 0 & 0 & 0 & 0\\
0 & 0 & 0 & 0 & 1 & 0 & 0 & 0 & 0 & 0\\
0 & 1 & 0 & 0 & 0 & 0 & 0 & 0 & 0 & 0\\
0 & 0 & 0 & 0 & 0 & 0 & 0 & 0 & 0 & 1\\
0 & 0 & 0 & 0 & 0 & 0 & 0 & 1 & 0 & 0\\
0 & 0 & 0 & 0 & 0 & 0 & 0 & 0 & 1 & 0\\
\end{array} \right].\]
\end{center}
Then $x_{n}$ and $y_{n}$ are elements of $G (= SL_{n}(q))$ of orders $2$ and $3$, respectively. Denote  
\begin{center}
\[z_{9} = x_{9}y_{9} = \left[ \begin{array}{ccccccccc}
0 & 0 & -1 & 0 & 0 & 0 & 0 & \alpha_{5} & \alpha_{5}\alpha_{8}^{-1}\\
-1 & 0 & 0 & 0 & 0 & 0 & 0 & \alpha_{4} & \alpha_{4}\alpha_{8}^{-1}\\
0 & 0 & 0 & 0 & 0 & -1 & 0 & \alpha_{6} & \alpha_{3}\alpha_{8}^{-1}\\
0 & -1 & 0 & 0 & 0 & 0 & 0 & \alpha_{3} & \alpha_{6}\alpha_{8}^{-1}\\
0 & 0 & 0 & -1 & 0 & 0 & 0 & \alpha_{2} & \alpha_{2}\alpha_{8}^{-1}\\
0 & 0 & 0 & 0 & 0 & 0 & -1 & \alpha_{7} & \alpha_{1}\alpha_{8}^{-1}\\
0 & 0 & 0 & 0 & 0 & 0 & 0 &  \alpha_{8} & 0\\
0 & 0 & 0 & 0 & -1 & 0 & 0 & \alpha_{1} & \alpha_{7}\alpha_{8}^{-1}\\
0 & 0 & 0 & 0 & 0 & 0 & 0 & 0 & \alpha_{8}^{-1}\\
 \end{array} \right],\]
\[z_{10} = x_{10}y_{10} = \left[ \begin{array}{cccccccccc}
0 & 0 & 0 & 0 & 0 & -1 & 0 & 0 & \alpha_{3} & \alpha_{2}\alpha_{9}^{-1}\\
0 & 0 & 0 & 0 & -1 & 0 & 0 & 0 & \alpha_{7} & \alpha_{4}\alpha_{9}^{-1}\\
0 & 0 & 0 & 0 & 0 & 0 & -1 & 0 & \alpha_{5} & \alpha_{5}\alpha_{9}^{-1}\\
-1 & 0 & 0 & 0 & 0 & 0 & 0 & 0 & \alpha_{2} & \alpha_{3}\alpha_{9}^{-1}\\
0 & 0 & 0 & 0 & 0 & 0 & 0 & -1 & \alpha_{8} & \alpha_{1}\alpha_{9}^{-1}\\
0 & 0 & -1 & 0 & 0 & 0 & 0 & 0 & \alpha_{4} & \alpha_{7}\alpha_{9}^{-1}\\
0 & -1 & 0 & 0 & 0 & 0 & 0 & 0 & \alpha_{6} & \alpha_{6}\alpha_{9}^{-1}\\
0 & 0 & 0 & 0 & 0 & 0 & 0 & 0 & \alpha_{9} & 0\\
0 & 0 & 0 & -1 & 0 & 0 & 0 & 0 & \alpha_{1} & \alpha_{8}\alpha_{9}^{-1}\\
0 & 0 & 0 & 0 & 0 & 0 & 0 & 0 & 0 & \alpha_{9}^{-1}\\ \end{array} \right].\]
\end{center}
The characteristic polynomial of $z_{n}$ is $f_{z_{n}}(t) = (t - \alpha_{n-1}^{-1})f_{n}(t)$ and the characteristic roots $\alpha_{n-1}^{-1}$, $\omega$, $\omega^{q}$, $\omega^{q^{2}}$, $\omega^{q^{3}}$, . . . , $\omega^{q^{n-3}}$, and $\omega^{q^{n-2}}$ of $z_{n}$ are pairwise distinct. Then, in $GL_{n}(q^{n-1})$, $z_{n}$ is conjugate to the matrix diag $(\alpha_{n-1}^{-1}$, $\omega$, $\omega^{q}$, $\omega^{q^{2}}$, $\omega^{q^{3}}$, . . . , $\omega^{q^{n-3}}$, $\omega^{q^{n-2}})$ and hence $z_{n}$ is an element of $SL_{n}(q)$ of order $Q$.\\
\indent\indent Let $H_{n}$ is the subgroup of $G (= SL_{n}(q))$ generated by the above elements $x_{n}$ and $y_{n}$.\\

\indent\indent\textbf{Lemma  2.} \emph{The group $H_{n}$ can not stabilize one-dimensional subspaces or hyperplanes of the space $V$ or equivalently $H_{n}$ acts irreducible on $V$.}\\

\indent\indent P r o o f. Assume that $W$ is an $H_{n}$-invariant subspace of $V$ and $k$ = dim $W$, $k=1$ or $n-1$.\\
\indent\indent Let first $k = 1$ and $0 \neq w \in W$. Then $y_{n}(w) = \lambda w$ where $\lambda \in F$ and $\lambda^{3} = 1$. This yields
\begin{center}
$w = \mu_{1}(v_{1} + \lambda^{2} v_{2} + \lambda v_{3}) + \mu_{2}(v_{4} + \lambda^{2} v_{5} + \lambda v_{6}) + \mu_{3}(v_{7} + \lambda^{2} v_{8} + \lambda v_{9})$  $(\mu_{1}, \mu_{2}, \mu_{3} \in F)$ 
\end{center}
if $n = 9$, and
\begin{center}
$w = \mu_{1}^{'}v_{1} + \mu_{2}^{'}(v_{2} + \lambda v_{3}) + \mu_{3}^{'}(\lambda v_{4} + v_{5} + \lambda^{2} v_{6}) + \mu_{2}^{'} \lambda^{2} v_{7} + \mu_{4}^{'}(\lambda v_{8} + v_{9} + \lambda^{2} v_{10}) (\mu_{i}^{'} \in F)$ 
\end{center}
if $n = 10$. Moreover $\mu_{1}^{'} = 0$ if $\lambda \neq 1$.\\
Now $x_{n}(w) = \nu w$ where $\nu = \pm 1$. This yields consecutively $\mu_{3} \neq 0$ , $\mu_{4}^{'} \neq 0$, $\alpha_{n - 1} = \lambda^{2}\nu$, and (in case $n = 9$)
\begin{enumerate}[(1)]
  \item \rule{0pt}{0pt}\vspace*{-12pt}
    \begin{equation*}
     \lambda\nu\mu_{1} + \mu_{2} = (\lambda\nu\alpha_{3} + \lambda\alpha_{6})\mu_{3},
    \end{equation*}
\end{enumerate}
\vspace{0pt}
\begin{enumerate}[(2)]
  \item \rule{0pt}{0pt}\vspace*{-12pt}
    \begin{equation*}
     \mu_{2} = (\alpha_{1} - \lambda\nu + \nu\alpha_{7})\mu_{3},
    \end{equation*}
\end{enumerate}
\vspace{0pt}
\begin{enumerate}[(3)]
  \item \rule{0pt}{0pt}\vspace*{-12pt}
    \begin{equation*}
     (\nu + 1)(\lambda\mu_{2} - \alpha_{2}\mu_{3}) = 0,
    \end{equation*}
\end{enumerate}
\vspace{0pt}
\begin{enumerate}[(4)]
  \item \rule{0pt}{0pt}\vspace*{-12pt}
    \begin{equation*}
     (\nu + 1)(\mu_{1} - \lambda\alpha_{5}\mu_{3}) = 0,
    \end{equation*}
\end{enumerate}
\vspace{0pt}
\begin{enumerate}[(5)]
  \item \rule{0pt}{0pt}\vspace*{-12pt}
    \begin{equation*}
    (\nu + 1)(\mu_{1} - \lambda^{2}\alpha_{4}\mu_{3}) = 0;
    \end{equation*}
\end{enumerate}
in case $n = 10$ we obtain the following relations:\\
\begin{enumerate}[(1.)]
  \item \rule{0pt}{0pt}\vspace*{-12pt}
    \begin{equation*}
     \mu_{1}^{'} = -\nu\lambda\mu_{3}^{'} + (\lambda\alpha_{3}\alpha_{9}^{-1} + \lambda^{2}\alpha_{2})\mu_{4}^{'},
    \end{equation*}
\end{enumerate}
\vspace{0pt}
\begin{enumerate}[(2.)]
  \item \rule{0pt}{0pt}\vspace*{-12pt}
    \begin{equation*}
     \mu_{2}^{'} = -\nu\lambda^{2}\mu_{3}^{'} + (\lambda\alpha_{7}\alpha_{9}^{-1} + \lambda^{2}\alpha_{4})\mu_{4}^{'},
    \end{equation*}
\end{enumerate}
\vspace{0pt}
\begin{enumerate}[(3.)]
  \item \rule{0pt}{0pt}\vspace*{-12pt}
    \begin{equation*}
     \mu_{3}^{'} = (-\nu + \lambda_{2}\alpha_{1} + \lambda\alpha_{8}\alpha_{9}^{-1})\mu_{4}^{'},
    \end{equation*}
\end{enumerate}
\vspace{0pt}
\begin{enumerate}[(4.)]
  \item \rule{0pt}{0pt}\vspace*{-12pt}
    \begin{equation*}
     (\nu + 1)(\mu_{2}^{'} - \lambda\alpha_{5}\mu_{4}^{'}) = 0,
    \end{equation*}
\end{enumerate}
\vspace{0pt}
\begin{enumerate}[(5.)]
  \item \rule{0pt}{0pt}\vspace*{-12pt}
    \begin{equation*}
    (\nu + 1)(\lambda^{2}\alpha_{6} - \alpha_{5}) = 0.
    \end{equation*}
\end{enumerate}
In particular, we have $\alpha_{n-1}^{3} = \nu$ and $\alpha_{n-1}^{6} = 1$. This is impossible if $q = 5$ or $q > 7$ since then $\alpha_{n-1}$ has order $q - 1$.\\
Next, let us continue with the case $n = 9$. According to our assumption ($q \neq 2, 4$) only two possibilities left: $q = 3$ (and $\alpha_{8} = 1$), $q = 7$ (and $\alpha_{8}^{3} = 1 \neq \alpha_{8}$). So $\nu = 1$, $\alpha_{8} = \lambda^{2}$ and (1), (2), (3), (4), (5) produce $\alpha_{1} = \lambda^{2}\alpha_{2} - \alpha_{7} + \lambda$, $\alpha_{3} = \lambda\alpha_{2} + \lambda^{2}\alpha_{4} - \alpha_{6}$ and $\alpha_{5} = \lambda\alpha_{4}$. Now $f_{9}(-1) =  (1 + \lambda + \lambda^{2})(1 + \alpha_{2} + \alpha_{4}) = 0$ both for $q = 3$ and $q = 7$, an impossibility as $f_{9}(t)$ is irreducible over the field $F$.\\
Lastly, we treat the case $n = 10$, and as $q > 4$, that is $q = 7$ (and $\alpha_{9}^{3} = 1 \neq \alpha_{9}$). Thus $\nu = 1$ and $\alpha_{9} = \lambda^{2} \neq 1$. So 
$\lambda \neq 1$, $\mu_{1}^{'} = 0$ and from $(1.)$, $(2.)$, $(3.)$, $(4.)$, $(5.)$ we can extract that $\alpha_{1} = \lambda^{2}\alpha_{2} + \lambda^{2}\alpha_{3} - \alpha_{8} + \lambda$, $\alpha_{5} = -\lambda^{2}\alpha_{2} - \lambda^{2}\alpha_{3} + \lambda\alpha_{4} + \lambda\alpha_{7}$ and $\alpha_{6} = -\alpha_{2} -\alpha_{3}+ \lambda^{2}\alpha_{4} + \lambda^{2}\alpha_{7}$. Then $f_{10}(-1) =  -(1 + \lambda + \lambda^{2})(1 + \alpha_{4} + \alpha_{7}) = 0$, again an impossibility as $f_{10}(t)$ is irreducible over the field $F$.\\
\indent\indent Now let $k=n-1$. The subspace $U$ of $V$ which is generated by the vectors $v_{1}$, $v_{2}$, $v_{3}$, . . .  , $v_{n-1}$ is $\left\langle {z_{n}}\right\rangle$-invariant. If $W \neq U$ then $U \cap W$ is $\left\langle {z_{n}}\right\rangle$-invariant and dim $(U \cap W) = n-2$. This means that the characteristic polynomial of $z_{n}|_{U \cap W}$ has degree $n-2$ and must divide $f_{z_{n}}(t)$ which is impossible as $f_{n}(t)$ is irreducible over $F$. Thus $W = U$ but obviously $U$ is not $\left\langle {y_{n}}\right\rangle$-invariant, a contradiction.\\
\indent\indent The lemma is proved. (Note that the statement is false if $q = 2$ or $4$ in both cases, and additionally if $q = 3$ in case $n = 10$.) \hfill $\square$\\

\indent\indent Now, as $H_{n} = \left\langle{x_{n},y_{n}}\right\rangle$ acts irreducible on the space $V$ and it has an element of order $Q$, we conclude (by Lemma 1) that $H_{n}$ can not be contained in any maximal subgroup of $G (= SL_{n}(q))$. Thus $H_{n} = G$ and $G = \left\langle {x_{n},y_{n}}\right\rangle$ is a $(2,3)$-generated group. Obviously $\overline{x_{n}}$ and $\overline{y_{n}}$ are elements of respective orders $2$ and $3$ for the group $\overline{G} = PSL_{n}(q)$, and $\overline{G} = \left\langle {\overline{x_{n}},\overline{y_{n}}}\right\rangle$ is a $(2,3)$-generated group too.\\

\indent\indent \textbf{2.2.} Now we proceed to prove the $(2,3)$ - generation of the remaining groups $SL_{9}(2)$, $SL_{9}(4)$, $SL_{10}(2)$, $SL_{10}(3)$ and $SL_{10}(4)$. Below we provide elements $x_{n}^{(q)}$ and $y_{n}^{(q)}$ of orders $2$ and $3$, respectively, for each one of the groups $SL_{n}(q)$ in this list, and prove that $\left\langle {x_{n}^{(q)},y_{n}^{(q)}}\right\rangle = SL_{n}(q)$. In our consideration, counting the orders of some elements in the corresponding groups $\left\langle {x_{n}^{(q)},y_{n}^{(q)}}\right\rangle$, we rely on the great possibilities of Magma Computational Algebra System. We also use the orders of the maximal subgroups of the groups in the list above. (The maximal subgroups of these groups are classified in \cite{1}.)\\

\indent\indent Take the following two matrices of $SL_{9}(2)$:
\[ x_{9}^{(2)} = \left[ \begin{array}{ccccccccc}
1 & 0 & 0 & 0 & 0 & 0 & 0 & 0 & 0\\
0 & 0 & 1 & 0 & 0 & 0 & 0 & 0 & 0\\
0 & 1 & 0 & 0 & 0 & 0 & 0 & 0 & 0\\
0 & 0 & 0 & 0 & 1 & 0 & 0 & 0 & 0\\
0 & 0 & 0 & 1 & 0 & 0 & 0 & 0 & 0\\
0 & 0 & 0 & 0 & 0 & 1 & 1 & 0 & 1\\
0 & 0 & 0 & 0 & 0 & 1 & 0 & 1 & 1\\
0 & 0 & 0 & 0 & 0 & 0 & 1 & 1 & 1\\
0 & 0 & 0 & 0 & 0 & 1 & 1 & 1 & 0\\
\end{array} \right],
y_{9}^{(2)} = \left[ \begin{array}{ccccccccc}
0 & 1 & 0 & 0 & 0 & 0 & 0 & 0 & 0\\
1 & 1 & 0 & 0 & 0 & 0 & 0 & 0 & 0\\
0 & 0 & 0 & 1 & 0 & 0 & 0 & 0 & 0\\
0 & 0 & 1 & 1 & 0 & 0 & 0 & 0 & 0\\
0 & 0 & 0 & 0 & 0 & 1 & 0 & 0 & 0\\
0 & 0 & 0 & 0 & 1 & 1 & 0 & 0 & 0\\
0 & 0 & 0 & 0 & 0 & 0 & 0 & 0 & 1\\
0 & 0 & 0 & 0 & 0 & 0 & 1 & 0 & 0\\
0 & 0 & 0 & 0 & 0 & 0 & 0 & 1 & 0\\
\end{array} \right].\]
Then $x_{9}^{(2)}$ and $y_{9}^{(2)}$ are elements of respective orders $2$ and $3$ in the group $SL_{9}(2)$, and  $x_{9}^{(2)}y_{9}^{(2)}$ has order $73$; also $x_{9}^{(2)}y_{9}^{(2)}(x_{9}^{(2)}(y_{9}^{(2)})^{2})^{2}$ is an element in $\left\langle {x_{9}^{(2)},y_{9}^{(2)}}\right\rangle$ of order $3.127$. Since in $SL_{9}(2)$ there is no maximal subgroup of order divisible by $73.127$ it follows that $SL_{9}(2) = \left\langle {x_{9}^{(2)},y_{9}^{(2)}}\right\rangle$.\\

\indent\indent Now continue with the desired  matrices of $SL_{9}(4)$:
\[ x_{9}^{(4)} = \left[ \begin{array}{ccccccccc}
0 & 1 & 0 & 0 & 0 & 0 & 0 & 0 & 0\\
1 & 0 & 0 & 0 & 0 & 0 & 0 & 0 & 0\\
0 & 0 & 0 & 1 & 0 & 0 & 0 & 0 & 0\\
0 & 0 & 1 & 0 & 0 & 0 & 0 & 0 & 0\\
0 & 0 & 0 & 0 & 1 & 0 & 0 & 0 & 0\\
0 & 0 & 0 & 0 & 0 & 0 & 1 & 0 & 0\\
0 & 0 & 0 & 0 & 0 & 1 & 0 & 0 & 0\\
0 & 0 & 0 & 0 & 0 & 0 & 0 & 1 & \eta\\
0 & 0 & 0 & 0 & 0 & 0 & 0 & 0 & 1\\
\end{array} \right],
y_{9}^{(4)} = \left[ \begin{array}{ccccccccc}
1 & 0 & 0 & 0 & 0 & 0 & 0 & 0 & 0\\
0 & 0 & 1 & 0 & 0 & 0 & 0 & 0 & 0\\
0 & 1 & 1 & 0 & 0 & 0 & 0 & 0 & 0\\
0 & 0 & 0 & 0 & 0 & 1 & 0 & 0 & 0\\
0 & 0 & 0 & 1 & 0 & 0 & 0 & 0 & 0\\
0 & 0 & 0 & 0 & 1 & 0 & 0 & 0 & 0\\
0 & 0 & 0 & 0 & 0 & 0 & 1 & 1 & 1\\
0 & 0 & 0 & 0 & 0 & 0 & 1 & 1 & 0\\
0 & 0 & 0 & 0 & 0 & 0 & 0 & 1 & 0\\
\end{array} \right].\]
(Here $\eta$ is a generator of $GF(4)^{*}$.)\\
Besides that $x_{9}^{(4)}$ and $y_{9}^{(4)}$ have orders $2$ and $3$, respectively, $x_{9}^{(4)}y_{9}^{(4)}$ has order $3.5.43.127$ and in $\left\langle {x_{9}^{(4)},y_{9}^{(4)}}\right\rangle$ the order of the following element 
\begin{center} 
$(x_{9}^{(4)}(y_{9}^{(4)})^{2})^{2}(x_{9}^{(4)}y_{9}^{(4)})^{3}x_{9}^{(4)}(y_{9}^{(4)})^{2}(x_{9}^{(4)}y_{9}^{(4)})^{2}x_{9}^{(4)}(y_{9}^{(4)})^{2}(x_{9}^{(4)}y_{9}^{(4)})^{2}x_{9}^{(4)}(y_{9}^{(4)})^{2}x_{9}^{(4)}y_{9}^{(4)}$ 
\end{center}
is $3.7.19.73$. But no one maximal subgroup of $SL_{9}(4)$ has order divisible by $43.73$. Thus $ SL_{9}(4) = \left\langle {x_{9}^{(4)},y_{9}^{(4)}}\right\rangle$ is a $(2,3)$ - generated group too.\\

\indent\indent Next, we consider the appropriate pair of elements in the group $SL_{10}(2)$:
\[ x_{10}^{(2)} = \left[ \begin{array}{cccccccccc}
0 & 1 & 0 & 0 & 0 & 0 & 0 & 0 & 0 & 0\\
1 & 0 & 0 & 0 & 0 & 0 & 0 & 0 & 0 & 0\\
0 & 0 & 0 & 1 & 0 & 0 & 0 & 0 & 0 & 0\\
0 & 0 & 1 & 0 & 0 & 0 & 0 & 0 & 0 & 0\\
0 & 0 & 0 & 0 & 0 & 1 & 0 & 0 & 0 & 0\\
0 & 0 & 0 & 0 & 1 & 0 & 0 & 0 & 0 & 0\\
0 & 0 & 0 & 0 & 0 & 0 & 1 & 1 & 0 & 1\\
0 & 0 & 0 & 0 & 0 & 0 & 1 & 0 & 1 & 1\\
0 & 0 & 0 & 0 & 0 & 0 & 0 & 1 & 1 & 1\\
0 & 0 & 0 & 0 & 0 & 0 & 1 & 1 & 1 & 0\\
\end{array} \right],
y_{10}^{(2)} = \left[ \begin{array}{cccccccccc}
1 & 0 & 0 & 0 & 0 & 0 & 0 & 0 & 0 & 0\\
0 & 0 & 1 & 0 & 0 & 0 & 0 & 0 & 0 & 0\\
0 & 1 & 1 & 0 & 0 & 0 & 0 & 0 & 0 & 0\\
0 & 0 & 0 & 0 & 1 & 0 & 0 & 0 & 0 & 0\\
0 & 0 & 0 & 1 & 1 & 0 & 0 & 0 & 0 & 0\\
0 & 0 & 0 & 0 & 0 & 0 & 1 & 0 & 0 & 0\\
0 & 0 & 0 & 0 & 0 & 1 & 1 & 0 & 0 & 0\\
0 & 0 & 0 & 0 & 0 & 0 & 0 & 0 & 0 & 1\\
0 & 0 & 0 & 0 & 0 & 0 & 0 & 1 & 0 & 0\\
0 & 0 & 0 & 0 & 0 & 0 & 0 & 0 & 1 & 0\\
\end{array} \right].\]
Here we obtain that the order of $x_{10}^{(2)}y_{10}^{(2)}$ is $3.11.31$ and $x_{10}^{(2)}y_{10}^{(2)}(x_{10}^{(2)}(y_{10}^{(2)})^{2})^{2}$ has order $73$. But there is no maximal subgroup in $SL_{10}(2)$ of order divisible by $11.73$. So $SL_{10}(2) = \left\langle {x_{10}^{(2)},y_{10}^{(2)}}\right\rangle$.\\

\indent\indent Further, let us deal with the group $SL_{10}(3)$ and choose:
\[ x_{10}^{(3)} = \left[ \begin{array}{cccccccccc}
0 & 1 & 0 & 0 & 0 & 0 & 0 & 0 & 0 & 0\\
1 & 0 & 0 & 0 & 0 & 0 & 0 & 0 & 0 & 0\\
0 & 0 & 1 & 0 & 0 & 0 & 0 & 0 & 0 & 0\\
0 & 0 & 0 & 0 & 1 & 0 & 0 & 0 & 0 & 0\\
0 & 0 & 0 & 1 & 0 & 0 & 0 & 0 & 0 & 0\\
0 & 0 & 0 & 0 & 0 & 1 & 0 & 0 & 0 & 0\\
0 & 0 & 0 & 0 & 0 & 0 & 0 & 1 & 0 & 0\\
0 & 0 & 0 & 0 & 0 & 0 & 1 & 0 & 0 & 0\\
0 & 0 & 0 & 0 & 0 & 0 & 0 & 0 & -1 & 1\\
0 & 0 & 0 & 0 & 0 & 0 & 0 & 0 & 0 & 1\\
\end{array} \right],
y_{10}^{(3)} = \left[ \begin{array}{cccccccccc}
1 & 0 & 0 & 0 & 0 & 0 & 0 & 0 & 0 & 0\\
0 & 0 & 0 & 1 & 0 & 0 & 0 & 0 & 0 & 0\\
0 & 1 & 0 & 0 & 0 & 0 & 0 & 0 & 0 & 0\\
0 & 0 & 1 & 0 & 0 & 0 & 0 & 0 & 0 & 0\\
0 & 0 & 0 & 0 & 0 & 0 & 1 & 0 & 0 & 0\\
0 & 0 & 0 & 0 & 1 & 0 & 0 & 0 & 0 & 0\\
0 & 0 & 0 & 0 & 0 & 1 & 0 & 0 & 0 & 0\\
0 & 0 & 0 & 0 & 0 & 0 & 0 & 0 & 1 & 1\\
0 & 0 & 0 & 0 & 0 & 0 & 0 & 1 & 0 & -1\\
0 & 0 & 0 & 0 & 0 & 0 & 0 & 0 & 1 & 0\\
\end{array} \right].\]
The product of the last two matrices has order $11^{2}.61$ and the following element 
\begin{center}
$(x_{10}^{(3)}y_{10}^{(3)})^{2}x_{10}^{(3)}(y_{10}^{(3)})^{2}(x_{10}^{(3)}y_{10}^{(3)})^{2}x_{10}^{(3)}(y_{10}^{(3)})^{2}x_{10}^{(3)}y_{10}^{(3)}x_{10}^{(3)}(y_{10}^{(3)})^{2}x_{10}^{(3)}y_{10}^{(3)}$
\end{center}
is of order $2.13.757$. Checking the orders of the maximal subgroups of $SL_{10}(3)$ we can see that no one of them is a multiple of $61.757$ which means that $SL_{10}(3) = \left\langle {x_{10}^{(3)},y_{10}^{(3)}}\right\rangle$.\\

\indent\indent Lastly, we finish with the proof of the $(2,3)$ - generation of the group $SL_{10}(4)$ by taking its elements: 
\[ x_{10}^{(4)} = \left[ \begin{array}{cccccccccc}
1 & 0 & 0 & 0 & 0 & 0 & 0 & 0 & 0 & 0\\
0 & 0 & 1 & 0 & 0 & 0 & 0 & 0 & 0 & 0\\
0 & 1 & 0 & 0 & 0 & 0 & 0 & 0 & 0 & 0\\
0 & 0 & 0 & 0 & 1 & 0 & 0 & 0 & 0 & 0\\
0 & 0 & 0 & 1 & 0 & 0 & 0 & 0 & 0 & 0\\
0 & 0 & 0 & 0 & 0 & 1 & 0 & 0 & 0 & 0\\
0 & 0 & 0 & 0 & 0 & 0 & 0 & 1 & 0 & 0\\
0 & 0 & 0 & 0 & 0 & 0 & 1 & 0 & 0 & 0\\
0 & 0 & 0 & 0 & 0 & 0 & 0 & 0 & 1 & \eta\\
0 & 0 & 0 & 0 & 0 & 0 & 0 & 0 & 0 & 1\\
\end{array} \right],
y_{10}^{(4)} = \left[ \begin{array}{cccccccccc}
0 & 1 & 0 & 0 & 0 & 0 & 0 & 0 & 0 & 0\\
1 & 1 & 0 & 0 & 0 & 0 & 0 & 0 & 0 & 0\\
0 & 0 & 0 & 1 & 0 & 0 & 0 & 0 & 0 & 0\\
0 & 0 & 1 & 1 & 0 & 0 & 0 & 0 & 0 & 0\\
0 & 0 & 0 & 0 & 0 & 0 & 1 & 0 & 0 & 0\\
0 & 0 & 0 & 0 & 1 & 0 & 0 & 0 & 0 & 0\\
0 & 0 & 0 & 0 & 0 & 1 & 0 & 0 & 0 & 0\\
0 & 0 & 0 & 0 & 0 & 0 & 0 & 1 & 1 & 1\\
0 & 0 & 0 & 0 & 0 & 0 & 0 & 1 & 1 & 0\\
0 & 0 & 0 & 0 & 0 & 0 & 0 & 0 & 1 & 0\\
\end{array} \right].\]
(Recall that $\eta$ is a generator of $GF(4)^{*}$.)\\
In this case $x_{10}^{(4)}y_{10}^{(4)}$ has order $3.19.73$ and the element $(x_{10}^{(4)}(y_{10}^{(4)})^{2})^{3}x_{10}^{(4)}y_{10}^{(4)}(x_{10}^{(4)}(y_{10}^{(4)})^{2})^{6}$ has order $5.11.31.41$. Similarly, as in the previous cases, we can conclude that $SL_{10}(4)$ is generated by the above matrices because no one of its maximal subgroups has order divisible by $41.73$.\\
\indent\indent Finally, it is obvious that the projective images (if it is necessary) of all these elements $x_{n}^{(q)}$ and $y_{n}^{(q)}$ (again of orders $2$ and $3$, respectively) generate the correspondent simple special linear group.\\ 

\indent\indent \textbf{Acknowledgement}. We express our gratitude to Prof. Marco Antonio Pellegrini who provided us with all these generators for the groups $SL_{9}(2)$, $SL_{9}(4)$, $SL_{10}(2)$, $SL_{10}(3)$ and $SL_{10}(4)$.\\

\indent\indent \textbf{2.3.} Finally let $G = SL_{11}(q)$ and $\overline{G} = G/Z(G) = PSL_{11}(q)$, where $q = p^{e}$ and $p$ is a prime. Set $d = (11,q - 1)$ and $Q=(q^{11}-1)/(q-1)$. It is easily seen that here $(6,Q)=1$. The group $G$ acts (naturally) on an eleven-dimensional vector space $V = F^{11}$ over the field $F = GF(q)$.\\
\indent\indent\ We shall make use of the known list of maximal subgroups of ${G}$ given in \cite{1}. In Aschbaher's notation any maximal subgroup of ${G}$ belongs to one of the following families \emph{$C_{1}, C_{2}, C_{3}, C_{5}, C_{6}, C_{8}$}, and \emph{S}. Roughly speaking, they are:
\begin{itemize}
\item \emph{$C_{1}$}: stabilizers of subspaces of $V$,
\item \emph{$C_{2}$}: stabilizers of direct sum decompositions of $V$,
\item \emph{$C_{3}$}: stabilizers of extension fields of $F$ of prime degree,
\item \emph{$C_{5}$}: stabilizers of subfields of $F$ of prime index,
\item \emph{$C_{6}$}: normalizers of extraspecial groups in absolutely irreducible representations,
\item \emph{$C_{8}$}: classical groups on $V$ contained in $G$,
\item \emph{S}: almost simple groups, absolutely irreducible on $V$, and the representation of their (simple) \emph{socles} on $V$ can not be realized over proper subfields of $F$; not continued in members of \emph{$C_{8}$}. 
\end{itemize}
In \cite{1} the representatives of the conjugacy classes of maximal subgroups of ${G}$ are specified in Tables $8.70$ and $8.71$. For the reader's convenience we provide the exact list of maximal subgroups of $G$ together with their orders. The notation used here for group structures is standard group-theoretic notation as in \cite{1}. Especially, $A \times B$ is the direct product of groups $A$ and $B$, and we write $A:B$ or $A.B$ to denote a split extension of $A$ by $B$ or an extension of $A$ by $B$ of unspecified type, respectively; the cyclic group of order $n$ is simple denoted by $n$, and $E_{q^{k}}$ stands for an elementary abelian group of order $q^{k}$.\\
If $M$ is a maximal subgroup of $G$ then one of the following holds.
\begin{enumerate}
\item $M \cong E_{q^{10}}:GL_{10}(q)$ of order $q^{55}(q - 1)(q^{2} - 1)(q^{3} - 1)(q^{4} - 1)(q^{5} - 1)(q^{6} - 1)(q^{7} - 1)(q^{8} - 1)(q^{9} - 1)(q^{10} - 1)$.
\item $M \cong E_{q^{18}}:(SL_{9}(q)\times SL_{2}(q)):(q - 1)$ of order $q^{55}(q - 1)(q^{2} - 1)^{2}(q^{3} - 1)(q^{4} - 1)(q^{5} - 1)(q^{6} - 1)(q^{7} - 1)(q^{8} - 1)(q^{9} - 1)$.
\item $M \cong E_{q^{24}}:(SL_{8}(q)\times SL_{3}(q)):(q - 1)$ of order $q^{55}(q - 1)(q^{2} - 1)^{2}(q^{3} - 1)^{2}(q^{4} - 1)(q^{5} - 1)(q^{6} - 1)(q^{7} - 1)(q^{8} - 1)$.
\item $M \cong E_{q^{28}}:(SL_{7}(q)\times SL_{4}(q)):(q - 1)$ of order $q^{55}(q - 1)(q^{2} - 1)^{2}(q^{3} - 1)^{2}(q^{4} - 1)^{2}(q^{5} - 1)(q^{6} - 1)(q^{7} - 1)$.
\item $M \cong E_{q^{30}}:(SL_{6}(q)\times SL_{5}(q)):(q - 1)$ of order $q^{55}(q - 1)(q^{2} - 1)^{2}(q^{3} - 1)^{2}(q^{4} - 1)^{2}(q^{5} - 1)^{2}(q^{6} - 1)$.
\item $M \cong (q - 1)^{10}: S_{11}$ (if $q \geq 5$) of order $2^{8}.3^{4}.5^{2}.7.11.(q - 1)^{10}$.
\item $M \cong \frac{q^{11} - 1}{q - 1}:11$ of order $11.\frac{q^{11} - 1}{q - 1}$.
\item $M \cong SL_{11}(q_{0}). (11,\frac{q - 1}{q_{0} - 1})$ (if $q = q_{0}^{r}$, $r$ prime) of order  $q_{0}^{55}(q_{0}^{2} - 1)(q_{0}^{3} - 1)(q_{0}^{4} - 1)(q_{0}^{5} - 1)(q_{0}^{6} - 1)(q_{0}^{7} - 1)(q_{0}^{8} - 1)(q_{0}^{9} - 1)(q_{0}^{10} - 1)(q_{0}^{11} - 1).(11,\frac{q - 1}{q_{0} - 1})$.
\item $M \cong 11_{+}^{1+2}:Sp_{2}(11)$ (if $q = p \equiv 1$ (mod $11$) or $q=p^{5}$ and $p \equiv 3, 4, 5, 9$ (mod $11$)) of order $2^{3}.3.5.11^{4}$ (here $11_{+}^{1+2}$ stands for an extraspecial group of order $11^{3}$ and exponent $11$).
\item $M \cong d \times SO_{11}(q)$ (if $q$ is odd) of order $d.q^{25}(q^{2} - 1)(q^{4} - 1)(q^{6} - 1)(q^{8} - 1)(q^{10} - 1)$.
\item $M \cong (11,q_{0} - 1) \times SU_{11}(q_{0})$ (if $q = q_{0}^{2}$) of order  $q_{0}^{55}(q_{0}^{2} - 1)(q_{0}^{3} + 1)(q_{0}^{4} - 1)(q_{0}^{5} + 1)(q_{0}^{6} - 1)(q_{0}^{7} + 1)(q_{0}^{8} - 1)(q_{0}^{9} + 1)(q_{0}^{10} - 1)(q_{0}^{11} + 1).(11,q_{0} - 1)$.
\item $M \cong d \times L_{2}(23)$ (if $q = p \equiv 1, 2, 3, 4, 6, 8, 9, 12, 13, 16, 18$ (mod $23$), $q \neq 2$) of order $2^{3}.3.11.23.d$.
\item $M \cong d \times U_{5}(2)$ (if $q = p \equiv 1$ (mod $3$)) of order  $2^{10}.3^{5}.5.11.d$.
\item $M \cong M_{24}$ (if $q = 2$) of order $2^{10}.3^{3}.5.7.11.23$.
\end{enumerate}

\indent\indent Now we proceed as follows. First we prove that there is only one type of maximal subgroups of $G$ whose order is a multiple of $Q$; actually these are the groups (in Aschbaher's class \emph{$C_{3}$}) in case $7$ above, of order $11.\frac{q^{11} - 1}{q - 1}$. Into the second step we find out two elements $x_{11}$ and $y_{11}$ of respective orders $2$ and $3$ in $G$ such that their product has got order $Q$. Finally, we deduce that the group $G$ is generated by these two elements. Then the projective images of these elements will generate the group $\overline{G}$.\\
\indent\indent Let us start with the first step in our strategy. In order to prove the above mentioned arithmetic fact we use the well-known Zsigmondy's theorem, and take a primitive prime divisor of $p^{11e} - 1$, i.e., a prime $r$ which divides $p^{11e} - 1$ but does not divide $p^{i} - 1$ for $0 < i < 11e$. Obviously $r \geq 23$ (as $r - 1$ is a multiple of $11e$) and also $r$ divides $Q$. Now it is easy to be seen that the only maximal subgroups of orders divisible by $r$ are those in cases $11$, and $12$ or $14$ with $r = 23$. In case $11$ if $Q = \frac{q_{0}^{22} - 1}{q_{0}^{2} - 1}$ divides the order of $M$, then $\frac{q_{0}^{11} - 1}{q_{0} - 1}$ should be a factor of the integer $q_{0}^{55}(q_{0} + 1)(q_{0}^{2} - 1)(q_{0}^{3} + 1)(q_{0}^{4} - 1)(q_{0}^{5} + 1)(q_{0}^{6} - 1)(q_{0}^{7} + 1)(q_{0}^{8} - 1)(q_{0}^{9} + 1)(q_{0}^{10} - 1).(11,q_{0} - 1)$, an impossibility, by the same Zsigmondy's theorem. As for the groups in case $12$, we have $Q = \frac{p^{11} - 1}{p - 1} \geq  \frac{3^{11} - 1}{3 - 1} > 2^{3}.3.11^{2}.23 \geq |M|$. Lastly, in case $14$ $Q = 2^{11} - 1 = 23.89$ does not divide the order of $M_{24}$.\\
\indent\indent Further, let us choose for $x_{11}$ the matrix
\begin{center}
\[ x_{11} = \left[ \begin{array}{rrrrrrrrrrr}
0 & 0 & 0 & 0 & 0 & 0 & 0 & 0 & 0 & 0 & 1\\
0 & 0 & 0 & 0 & 0 & 0 & 0 & 0 & 0 & 1 & 0\\
0 & 0 & 0 & 0 & 0 & 0 & 0 & 0 & 1 & 0 & 0\\
0 & 0 & 0 & 0 & 0 & 0 & 0 & 1 & 0 & 0 & 0\\
0 & 0 & 0 & 0 & 0 & 0 & 1 & 0 & 0 & 0 & 0\\
0 & 0 & 0 & 0 & 0 & -1 & 0 & 0 & 0 & 0 & 0\\
0 & 0 & 0 & 0 & 1 & 0 & 0 & 0 & 0 & 0 & 0\\
0 & 0 & 0 & 1 & 0 & 0 & 0 & 0 & 0 & 0 & 0\\
0 & 0 & 1 & 0 & 0 & 0 & 0 & 0 & 0 & 0 & 0\\
0 & 1 & 0 & 0 & 0 & 0 & 0 & 0 & 0 & 0 & 0\\
1 & 0 & 0 & 0 & 0 & 0 & 0 & 0 & 0 & 0 & 0\\
 \end{array} \right]\] 
\end{center}
and $y_{11}$ to be in the form
\begin{center}
\[ y_{11} = \left[ \begin{array}{rrrrrrrrrrl}
-1 & -1 & 0 & 0 & 0 & 0 & 0 & 0 & 0 & 0 & \delta_{1}\\
1 & 0 & 0 & 0 & 0 & 0 & 0 & 0 & 0 & 0 & \delta_{2}\\
0 & 0 & -1 & -1 & 0 & 0 & 0 & 0 & 0 & 0 & \delta_{3}\\
0 & 0 & 1 & 0 & 0 & 0 & 0 & 0 & 0 & 0 & \delta_{4}\\
0 & 0 & 0 & 0 & -1 & -1 & 0 & 0 & 0 & 0 & \delta_{5}\\
0 & 0 & 0 & 0 & 1 & 0 & 0 & 0 & 0 & 0 & \delta_{6}\\
0 & 0 & 0 & 0 & 0 & 0 & -1 & -1 & 0 & 0 & \delta_{7}\\
0 & 0 & 0 & 0 & 0 & 0 & 1 & 0 & 0 & 0 & \delta_{8}\\
0 & 0 & 0 & 0 & 0 & 0 & 0 & 0 & -1 & -1 & \delta_{9}\\
0 & 0 & 0 & 0 & 0 & 0 & 0 & 0 & 1 & 0 & \delta_{10}\\
0 & 0 & 0 & 0 & 0 & 0 & 0 & 0 & 0 & 0 & 1\\
 \end{array} \right].\] 
\end{center}
Then $x_{11}$ is an involution of $G$ and $y_{11}$ is an element of order $3$ in $G$ for any $\delta_{1}$, $\delta_{2}$, $\delta_{3}$, $\delta_{4}$, $\delta_{5}$, $\delta_{6}$, $\delta_{7}$, $\delta_{8}$, $\delta_{9}$, $\delta_{10} \in GF(q)$, also 
\begin{center}
\[ z_{11}=x_{11}y_{11} = \left[ \begin{array}{rrrrrrrrrrr}
0 & 0 & 0 & 0 & 0 & 0 & 0 & 0 & 0 & 0 & 1\\
0 & 0 & 0 & 0 & 0 & 0 & 0 & 0 & 1 & 0 & \delta_{10}\\
0 & 0 & 0 & 0 & 0 & 0 & 0 & 0 & -1 & -1 & \delta_{9}\\
0 & 0 & 0 & 0 & 0 & 0 & 1 & 0 & 0 & 0 & \delta_{8}\\
0 & 0 & 0 & 0 & 0 & 0 & -1 & -1 & 0 & 0 & \delta_{7}\\
0 & 0 & 0 & 0 & -1 & 0 & 0 & 0 & 0 & 0 & -\delta_{6}\\
0 & 0 & 0 & 0 & -1 & -1 & 0 & 0 & 0 & 0 & \delta_{5}\\
0 & 0 & 1 & 0 & 0 & 0 & 0 & 0 & 0 & 0 & \delta_{4}\\
0 & 0 & -1 & -1 & 0 & 0 & 0 & 0 & 0 & 0 & \delta_{3}\\
1 & 0 & 0 & 0 & 0 & 0 & 0 & 0 & 0 & 0 & \delta_{2}\\
-1 & -1 & 0 & 0 & 0 & 0 & 0 & 0 & 0 & 0 & \delta_{1}\\
 \end{array} \right].\] 
\end{center}
The characteristic polynomial of $z_{11}$ is
\begin{center}
$f_{z_{11}}(t) = t^{11}-\delta_{1}t^{10}+(\delta_{10}-1)t^{9}+(2\delta_{1}+\delta_{3}+1)t^{8}-(\delta_{1}+\delta_{8}+\delta_{9}+2\delta_{10}+1)t^{7}-(\delta_{1}-\delta_{2}+\delta_{3}+\delta_{5}-\delta_{10})t^{6}+(\delta_{1}+\delta_{3}-\delta_{6}+\delta_{7}+\delta_{8}+\delta_{9}+\delta_{10}+1)t^{5}+(\delta_{1}-\delta_{2}-\delta_{4}-\delta_{7}-\delta_{8}-\delta_{9}-\delta_{10})t^{4}-(\delta_{1}-\delta_{2}-\delta_{4}+\delta_{9}+\delta_{10}+2)t^{3}+(\delta_{2}+\delta_{9}+\delta_{10}+2)t^{2}-(\delta_{2}-1)t-1$
\end{center}
Now let take an element $\omega$  of order $Q$ in the multiplicative group of the field $GF(q^{11})$ and put 
\begin{center}
$l(t) = (t - \omega)(t - \omega^{q})(t - \omega^{q^{2}})(t - \omega^{q^{3}})(t - \omega^{q^{4}})(t - \omega^{q^{5}})(t - \omega^{q^{6}})(t - \omega^{q^{7}})(t - \omega^{q^{8}})(t - \omega^{q^{9}})(t - \omega^{q^{10}}) = t^{11} - at^{10} + bt^{9} - ct^{8} + dt^{7} - et^{6} + ft^{5}-gt^{4}+ht^{3}-kt^{2}+mt- 1$.
\end{center}
The last polynomial has all its coefficients in the field $GF(q)$ and the roots of $l(t)$ are pairwise distinct (in fact, the polynomial $l(t)$ is irreducible over $GF(q)$ which is not necessary for our considerations). The polynomials $f_{z_{11}}(t)$ and $l(t)$ are identically equal if  
\begin{center}
$\delta_{1}=a$, $\delta_{2}=-m+1$, $\delta_{3}=-2a-c-1$, $\delta_{4}=a+2m-2-k+h$, $\delta_{5}=a-m+3+c+b+e$, $\delta_{6}=-a-c+1+g-m+k-h-f$, $\delta_{7}=3-m+k-h+g+a+b+d$, $\delta_{8}=-a+k-m+1-b-d$, $\delta_{9}=m-4-b-k$, $\delta_{10}= b+1$
\end{center}
For these values of $\delta_{i} (i=1,...,10)$ $f_{z_{11}}(t)=l(t)$ and then, in $GL_{11}(q^{11})$, $z_{11}$ is conjugate to diag $(\omega, \omega^{q}, \omega^{q^{2}}, \omega^{q^{3}}, \omega^{q^{4}}, \omega^{q^{5}}, \omega^{q^{6}}, \omega^{q^{7}}, \omega^{q^{8}}, \omega^{q^{9}}, \omega^{q^{10}})$ and hence $z_{11}$ is an element of $G$ of order $Q$.\\
\indent\indent Now, $H = \left\langle x_{11},y_{11}\right\rangle$ is a subgroup of $G$ of order divisible by $6Q$. We have already proved above that the only maximal subgroup of $G$ whose order is a multiple of $Q$ is that in Aschbaher's class \emph{$C_{3}$}, of order $11Q$, which means that $H$ can not be contained in any maximal subgroup of $G$. Thus $H = G$ and $G = \left\langle x_{11},y_{11}\right\rangle$ is a $(2,3)$-generated group; $\overline{G} = \left\langle {\overline{x_{11}},\overline{y_{11}}}\right\rangle$ is a $(2,3)$-generated group too. \\
\indent\indent This completes the proof of the theorem.\hfill $\square$\\
\begin{center}

\end{center}

\vspace{25pt}

E. Gencheva and Ts. Genchev\\
Department of Mathematics\\
Technical University of Varna \\
Varna, Bulgaria\\
e-mail: elenkag@abv.bg; genchev57@yahoo.com\\

\vspace{40pt}

K. Tabakov\\
Faculty of Mathematics and Informatics\\
Department of Algebra\\
"St. Kliment Ohridski" University of Sofia\\
Sofia, Bulgaria\\
e-mail: ktabakov@fmi.uni-sofia.bg\\
\end{document}